\documentclass[english]{article}
\usepackage[T1]{fontenc}
\usepackage[latin9]{inputenc}
\usepackage{geometry}
\usepackage{amsmath}
\usepackage{amsthm}
\usepackage{amssymb}

\makeatletter
\newcommand{\lyxaddress}[1]{
\par {\raggedright #1
\vspace{1.4em}
\noindent\par}
}
\theoremstyle{plain}
\newtheorem{thm}{\protect\theoremname}
\theoremstyle{remark}
\newtheorem{rem}[thm]{\protect\remarkname}
\theoremstyle{plain}
\newtheorem{cor}[thm]{\protect\corollaryname}

\makeatother

\usepackage{babel}
\providecommand{\corollaryname}{Corollary}
\providecommand{\remarkname}{Remark}
\providecommand{\theoremname}{Theorem}

\begin{document}

\title{On the extraction of the boundary conditions and the boundary ports
in second-order field theories}

\author{Markus Schöberl and Kurt Schlacher}
\maketitle

\lyxaddress{Institute of Automatic Control and Control Systems Technology, Johannes
Kepler University, Altenbergerstr. 69, 4040 Linz, Austria}
\begin{abstract}
In this paper we consider second-order field theories in a variational
setting. From the variational principle the Euler-Lagrange equations
follow in an unambiguous way, but it is well known that this is not
true for the Cartan form. This has also consequences on the derivation
of the boundary conditions when non trivial variations are allowed
on the boundary. By posing extra conditions on the set of possible
boundary terms we exploit the degree of freedom in the Cartan form
to extract physical meaningful boundary expressions. The same mathematical
machinery will be applied to derive the boundary ports in a Hamiltonian
representation of the partial differential equations which is crucial
for energy based control approaches. Our results will be visualized
for mechanical systems such as beam and plate models.
\end{abstract}

\section{Introduction}

Many distributed-parameter systems described by partial differential
equations (pdes) allow for a variational formulation such that they
can be analyzed based on the corresponding Lagrangian density. In
a geometric setting the Euler-Lagrange equations are constructed by
using the so-called Cartan form, which provides a coordinate free
way to perform integration by parts, see e.g. \cite{Aldaya,GotayExt,Krupka,SaundersCrampin,Saunders}
and references therein. It is well-known that a unique formulation
for the Euler-Lagrange equations can be found but this is not true
for the Cartan form when the variational problem involves derivative
variables of second-order or higher and more than one independent
coordinate. This is connected to the commutativity of repeated partial
differentiation. In the literature many attempts have been made to
fix the Cartan form also in the higher-order field theories, e.g.
by using connections, see e.g \cite{Garcia} or to avoid the ambiguities
by using different approaches e.g. the Skinner and Rusk formalism,
see e.g. \cite{Campos}. Another consequence of the ambiguities of
the Cartan form (i.e. the need for repeated integration by parts)
is the fact that also the terms in the boundary integral are affected
by this non-uniqueness, when non-trivial variations on the boundary
are explicitly allowed. In the mathematical physicist literature this
effect is usually ignored and only trivial boundary conditions are
considered, but e.g. for mechanical structures such as beams, plates
or shells it is of great importance to allow for non trivial boundary
terms and consequently also for non-zero energy flow over the boundary,
especially for engineering applications, e.g. for control purposes
and especially for boundary control systems.

In this paper we consider non-trivial variations on the boundary and
we restrict ourselves to second-order field theories. The case of
second-order field theories is in some sense special, see e.g. \cite{KouranbaevaShkoller,PrietoMartinezRomanRoj,SaundersCrampin},
but to our best knowledge a systematic derivation of non-trivial boundary
conditions based on the Cartan form has not been reported in the literature
so far. In \cite{KouranbaevaShkoller} it is shown by evaluating the
variational principle (using integration by parts) that the boundary
term is indeed connected with the Cartan form but the extraction of
the meaningful boundary conditions is not discussed in this setting.
In \cite{Moreno} natural boundary conditions for higher-order field
theories in a variational setting are discussed based on flag fibrations
and again using integration by parts.

We demonstrate that by posing extra conditions on the set of possible
boundary terms, one is able to derive the physically admissible boundary
conditions which will be demonstrated for mechanical systems where
the Lagrangian depends on second-order derivatives. To this end, we
use the ambiguities in constructing the Cartan form and we exploit
the special structure of the boundary manifold of the underlying variational
problem, where we make use of an adapted coordinate system as it has
been suggested e.g. in \cite{EnnsPhD,Moreno} to present our result
in a compact and simple manner. Some preliminary work has been reported
in \cite{SchoeberlMathmod15} and \cite{SchoeberlLHMNLC}.

Furthermore, we will focus on partial differential equations (pdes)
that allow for a Hamiltonian formulation. Given pdes in a Lagrangian
setting, different approaches to construct Hamiltonian counterparts
have been discussed in the literature including symplectic, polysymplectic
and multisymplectic cases see, e.g. \cite{Giachetta,Gotay,Olver}.
From a control and a system theoretic point of view the class of port-Hamiltonian
systems with dissipation \cite{vanderSchaft} is a very prominent
system class, which has been extended from the lumped-parameter case
(systems described by ordinary differential equations) to the distributed-parameter
pde case. However, different port-Hamiltonian formulations exist,
see e.g. \cite{SchaftMaschke,Maschke2005,Machelli2004I,Machelli2004II}
for an approach based on a Stokes-Dirac structure whereas in \cite{Enns,SchlacherHam2008,SchoeberlLHMNLC,SchoeberlTAC2013,SchoeberlAut2014PDE}
an alternative strategy based on jet-bundles has been considered.
In this paper we will consider the jet-bundle approach and apply the
theory of Cartan forms to be able to extract the correct boundary
ports when the Hamiltonian depends on second-order derivative variables
in a formal and systematic way.

The main contribution of the paper is i) to provide a simple and algorithmic
procedure to derive the physical admissible boundary conditions for
variational problems by means of the Cartan form where the Lagrangian
depends on second-order derivative variables and ii) to apply this
mathematical machinery developed for the Lagrangian scenario to a
port-Hamiltonian representation in order to analyze the power balance
relation. The paper is organized as follows. In Section 2 some geometric
fundamentals are presented and in the third section we discuss first-order
theories where we consider the well-known case of the Cartan form
approach to derive the pdes and the boundary terms in a Lagrangian
setting. Then, Section 4 is devoted to the second-order case in the
Lagrangian framework, where we show how to extract the boundary conditions
by using the presented machinery. In Section 5 we discuss the Euler-Bernoulli
beam and the Kirchhoff plate in the Lagrangian framework. In the sixth
section the port-Hamiltonian perspective is analyzed within this mathematical
setting and the Kirchhoff plate example is rearranged in a Hamiltonian
formulation to make the power ports visible. 

\section{Notation and Preliminaries}

In this paper we will apply differential geometric methods, in particular
we use jet-bundles and induced geometric objects, see \cite{Saunders}
for an extensive treatise. We use tensor notation and especially Einsteins
convention on sums where we will not indicate the range of the indices
when they are clear from the context. Given a manifold $\mathcal{M}$
with local coordinates $(x^{\alpha})$ where $\alpha=1,\ldots,\mathrm{dim}(\mathcal{M})$
we denote with $\mathcal{T}(\mathcal{M})\rightarrow\mathcal{M}$ and
$\mathcal{T}^{*}(\mathcal{M})\rightarrow\mathcal{M}$ the tangent
bundle and the cotangent bundle, respectively. Local coordinates for
$\mathcal{T}(\mathcal{M})$ and $\mathcal{T}^{*}(\mathcal{M})$ are
$(x^{\alpha},\dot{x}^{\alpha})$ as well as $(x^{\alpha},\dot{x}_{\alpha})$
with respect to the holonomic fibre basis $\{\partial_{\alpha}\}$
and $\{\mathrm{d}x^{\alpha}\}$. Furthermore, we use $\wedge$ for
the exterior product (wedge product), $\mathrm{d}$ is the exterior
derivative and $\rfloor$ denotes the natural contraction between
tensor fields. With $\mathrm{L}_{v}(\omega)$ we denote the Lie-derivative
of the geometric object $\omega$ (e.g. a differential form) with
respect to the vector field $v$. Furthermore, $C^{\infty}(\cdot)$
denotes the set of smooth functions on the corresponding manifold.

Let $\mathcal{X}\overset{\pi}{\rightarrow}\mathcal{D}$ be a bundle
with $\mathrm{dim}(\mathcal{X})=m+r$ where $\mathcal{D}$ is an orientable
$r-$dimensional manifold. Local coordinates for $\mathcal{X}$ adapted
to the fibration are $(X^{i},x^{\alpha})$ with $x^{\alpha}\,,\,\alpha=1,\ldots,m$
where $X^{i}\,,\,i=1,\ldots,r$ is a coordinate chart for $\mathcal{D}$.
We fix the volume form $\Omega$ on the manifold $\mathcal{D}$ as
$\Omega=\mathrm{d}X^{1}\wedge\ldots\wedge\mathrm{d}X^{r}$ and we
set $\Omega_{i}=\partial_{i}\rfloor\text{\ensuremath{\Omega}}$. Let
$\sigma$ be a section of the bundle $\mathcal{X}\rightarrow\mathcal{D}$,
i.e. a map $\sigma:\mathcal{D}\rightarrow\mathcal{X}$ with coordinate
representation $\sigma(X^{i})=(X^{i},\sigma^{\alpha}(X^{i}))$. The
set of sections of the bundle $\mathcal{X}\overset{\pi}{\rightarrow}\mathcal{D}$
will be denoted as $\mbox{\ensuremath{\Gamma}(\ensuremath{\pi})}$.
Then the $k-$th order partial derivatives of $\sigma^{\alpha}$ are
\[
\frac{\partial^{k}\sigma^{\alpha}}{(\partial X^{1})^{j_{1}}\cdots(\partial X^{r})^{j_{r}}}=\partial_{J}\sigma^{\alpha}
\]
where $J$ is an ordered multi-index with $\sum j_{i}=k=\#J$. The
components of the multi-index $J$ are $J(i)$ and $1_{k}$ is the
special multi index defined as $J(i)=\delta_{k}^{i}$ such that all
components of $1_{k}$ are zero except the $k-$th entry which is
one. Two multi-indices $J$ and $I$ can be added and subtracted component
wise, where the result of a subtraction may not be a multi-index.

Given two sections, such that at a point $p\in\mathcal{D}$ the sections
and their partial derivatives up to order $k$ coincide, these sections
are called $k-$equivalent. Given now a smooth local section $\sigma\in\mbox{\ensuremath{\Gamma_{p}}(\ensuremath{\pi})}$
at point $p\in\mathcal{D}$, the equivalence class of $n-$equivalent
sections that contains $\sigma$ is denoted $j_{p}^{n}(\sigma)$ (called
the $n-$jet at $p$) and the set of all these $n-$jets is the $n-$th
order jet-manifold $\mathcal{J}^{n}(\mathcal{X})$ with adapted coordinates
$(X^{i},x_{J}^{\alpha})$ with $0\leq\#J\leq n$. Furthermore, we
have $\pi_{k,l}:\mathcal{J}^{k}(\mathcal{X})\rightarrow\mathcal{J}^{l}(\mathcal{X})$
for $k>l$ as well as $\pi_{k,0}:\mathcal{J}^{k}(\mathcal{X})\rightarrow\mathcal{X}$
and $\pi_{k}:\mathcal{J}^{k}(\mathcal{X})\rightarrow\mathcal{D}$.
The tangent and the cotangent bundles of $\mathcal{J}^{k}(\mathcal{X})$
can be easily constructed and the holonomic fibre bases for $\mathcal{T}(\mathcal{J}^{k}(\mathcal{X}))$
and $\mathcal{T}^{*}(\mathcal{J}^{k}(\mathcal{X}))$ read as $\{\partial_{i},\,\partial_{\alpha}^{J}\}$
and $\{\mathrm{d}X^{i},\,\mathrm{d}x_{J}^{\alpha}\}\,$, where $0\leq\#J\leq k$,
respectively, with $\partial_{i}=\partial/\partial X^{i}$ and $\partial_{\alpha}^{J}=\partial/\partial x_{J}^{\alpha}$.

The total derivative $d_{1_{i}}$ with respect to the independent
variable $X^{i}$ reads as
\begin{equation}
d_{1_{i}}=\partial_{i}+x_{J+1_{i}}^{\alpha}\partial_{\alpha}^{J}\label{eq:TotalDer}
\end{equation}
and $(d_{1_{i}}f)\circ j^{n+1}(\sigma)=\partial_{i}(f\circ j^{n}(\sigma))$
is met by construction for all $f\in C^{\infty}(\mathcal{J}^{n}(\mathcal{X}))$
and $\sigma\in\Gamma(\pi)$. The $k-$th order total derivative meets
$d_{1_{i}}\in\pi_{k+1,k}^{*}\left(\mathcal{T}(\mathcal{J}^{k}(\mathcal{X}))\right)$.
The dual elements to the total derivatives are the so-called contact
forms. In local coordinates we have
\begin{equation}
\omega_{J}^{\alpha}=\mathrm{d}x_{J}^{\alpha}-x_{J+1_{i}}^{\alpha}\mathrm{d}X^{i}\,,\,\,\,\,\#J\geq0,\label{eq:ContactForm}
\end{equation}
where a $k-$th order contact form meets $\omega_{J}^{\alpha}\in\pi_{k+1,k}^{*}\left(\mathcal{T}^{*}(\mathcal{J}^{k}(\mathcal{X}))\right)$
with $\#J=k$ and $j^{k}(\sigma)^{*}(\omega_{J}^{\alpha})=0$ is met
for $\sigma\in\Gamma(\pi)$. Furthermore, based on the structure of
the underlying fibre bundle the operator of horizontalisation $\mathrm{hor}$,
see e.g. \cite{Krupka}, can be introduced as 
\[
\mathrm{hor}(\mathrm{d}X^{i})=\mathrm{d}X^{i},\,\mathrm{hor}(\mathrm{d}x_{J}^{\alpha})=x_{J+1_{i}}^{\alpha}\mathrm{d}X^{i}
\]
which has the property that $\mathrm{hor}(\omega_{J}^{\alpha})=0$. 

The geometric version of the calculus of variations rests on the use
of so-called variational vector fields $\eta$, whose flow is used
to perform variations of sections $\sigma\in\Gamma(\pi)$, where we
consider the special case of vertical variations. Hence, these vector
fields $\eta$ are tangent to the fibration of $\mathcal{X}\overset{\pi}{\rightarrow}\mathcal{D}$
and they are referred to as vertical vector fields, i.e. $\text{\ensuremath{\eta}}\in\mathcal{V}(\mathcal{X})$
with the coordinate representation $\eta=\eta^{\alpha}(X,x)\partial_{\alpha}$,
i.e. meeting $\dot{X}=0$. Furthermore, we need the $n-$th jet-prolongation
of a vertical vector field $\eta\in\mathcal{V}(\mathcal{X})$ which
is given as
\begin{equation}
j^{n}(\eta)=\eta^{\alpha}\partial_{\alpha}+d_{J}(\eta^{\alpha})\partial_{\alpha}^{J}\,,\,\,\,d_{J}=(d_{1_{1}})^{j_{1}}\ldots(d_{1_{r}})^{j_{r}}\label{eq:JetProv}
\end{equation}
with $1\leq\#J\leq n$, see e.g. \cite{Olver,Saunders}. It is well-known
that these jet-prolongations of vertical vector fields will play a
prominent role when formulating a variation problem as recapitulated
next, see again \cite{Saunders}.

Let us consider a variational problem on the bundle $\mathcal{X}\overset{\pi}{\rightarrow}\mathcal{D}$
where the Lagrangian is of the form $\mathfrak{L}=\mathcal{L}\Omega$
with $\mathcal{L}\in C^{\infty}(\mathcal{J}^{k}(\mathcal{X}))$ and
induces a functional $L(\sigma)=\int_{\mathcal{S}}j^{k}(\sigma)^{*}\mathfrak{L}$
for sections $\sigma\in\Gamma(\pi)$ where $\mathcal{S}$ is a submanifold
of $\mathcal{D}$. In the calculus of variations, the task is to find
so-called critical sections $\sigma$ that extremize the functional
$L(\sigma)$. To find these critical sections we will consider vertical
variations $\phi_{\epsilon}$ (fiber preserving maps described by
a flow). A critical section $\sigma$ meets
\[
\left.d_{\epsilon}(\int_{\mathcal{S}}(j^{k}(\phi_{\epsilon}\circ\sigma)^{*}\mathfrak{L})\right|_{\epsilon=0}=0,
\]
which is equivalent to
\begin{equation}
\int_{\mathcal{S}}j^{k}(\sigma)^{*}(\mathrm{L}_{j^{k}(\eta)}(\mathfrak{L}))=0,\label{eq:LagFirstOrderLie}
\end{equation}
see e.g. \cite{Saunders}, where $\eta$ is the infinitesimal generator
of the flow $\phi_{\epsilon}$. 

\section{First-order Case\label{sec:First-order-Case}}

In this section we will recapitulate the well-known facts of a variational
problem for first-order Lagrangians, i.e. $\mathcal{L}\in C^{\infty}(\mathcal{J}^{1}(\mathcal{X}))$.
However, contrary to most of the differential geometric approaches
in the literature we allow for non-trivial boundary conditions and
we will demonstrate how the Cartan form is connected to these boundary
terms. 
\begin{rem}
These non-trivial boundary terms play a crucial role in many applications,
e.g. in engineering, where systems can be interconnected via their
boundaries, see e.g. \cite{Enns,SchlacherHam2008,SchoeberlMCMDS2010,SchoeberlTAC2013,SchoeberlAut2014PDE}
for this port based approach, or where actuators directly influence
the boundary variables - boundary control systems.
\end{rem}
This introductory section serves mainly as a reminder but it should
clarify some of the ideas that are then used in second-order case
$k=2$ in the forthcoming section. Let us consider the equation (\ref{eq:LagFirstOrderLie})
for $k=1$ then using 
\[
\mathrm{L}_{j^{1}(\eta)}(\mathfrak{L})=j^{1}(\eta)\rfloor\mathrm{d}\mathfrak{L}+\mathrm{d}(j^{1}(\eta)\mathfrak{\rfloor L})
\]
there is a natural decomposition into a domain and a boundary term,
as by Stokes' theorem $\mathrm{d}(j^{1}(\eta)\mathfrak{\rfloor L})$
contributes to the boundary. Unfortunately, we cannot pick the domain
conditions directly as $j^{1}(\eta)\rfloor\mathrm{d}\mathfrak{L}$
may depend on derivatives of the variational field $\eta$ that induces
further terms on the boundary. However, if we are able to reformulate
the problem such that no prolongation of $\eta$ is necessary we are
done. The well-known solution of this problem is based on the use
of contact forms (\ref{eq:ContactForm}) and the Cartan form, see
e.g. \cite{Saunders}, with 
\begin{equation}
\mathfrak{C}=\mathcal{L}\Omega+\rho_{\alpha}^{j}\,\omega^{\alpha}\wedge\Omega_{j}\,,\,\,\,\omega^{\alpha}=\mathrm{d}x^{\alpha}-x_{1_{i}}^{\alpha}\mathrm{d}X^{i}\label{eq:Cartan1odansatz}
\end{equation}
implying $\int_{S}j^{1}(\sigma)^{*}\mathfrak{L}=\int_{\mathcal{S}}j^{1}(\sigma)^{*}\mathfrak{C}$
as $\mathrm{hor}(\mathfrak{C})=\mathfrak{L}$ is met, where the derivation
of the functions $\rho_{\alpha}^{j}$ will be described next. Instead
of (\ref{eq:LagFirstOrderLie}) we consider 
\[
\int_{\mathcal{S}}j^{1}(\sigma)^{*}(\mathrm{L}_{j^{1}(\eta)}(\mathfrak{C}))=0
\]
and $j^{1}(\eta)\rfloor\mathrm{d}\mathfrak{C}$ will deliver the domain
conditions (pdes) and from $j^{1}(\eta)\mathfrak{\rfloor C}$ the
boundary expressions can be extracted, if we choose the functions
$\rho_{\alpha}^{i}$ such that 
\[
\mathrm{hor}(j^{1}(\eta)\rfloor\mathrm{d}\mathfrak{C})=\mathrm{hor}(\eta\rfloor\mathrm{d}\mathfrak{C})
\]
is met. Consequently, 
\[
\mathrm{hor}(j^{1}(\eta)\rfloor\mathrm{d}\mathfrak{C})=\mathrm{hor}\left(j^{1}(\eta)\rfloor(\mathrm{d\mathcal{L}\wedge\Omega}+\mathrm{d}\rho_{\alpha}^{j}\wedge\omega^{\alpha}\wedge\Omega_{j}-\rho_{\alpha}^{j}\mathrm{d}x_{1_{j}}^{\alpha}\Omega)\right)
\]
follows, which can be evaluated as
\begin{equation}
\mathrm{hor}(j^{1}(\eta)\rfloor\mathrm{d}\mathfrak{C})=\mathrm{(\partial_{\alpha}\mathcal{L}}\eta^{\alpha}+\partial_{\alpha}^{1_{i}}\mathcal{L}\eta_{1_{i}}^{\alpha}-\rho_{\alpha}^{i}v_{1_{i}}^{\alpha})\Omega-\mathrm{hor}(\mathrm{d}\rho_{\alpha}^{i})\wedge\eta^{\alpha}\Omega_{i}\label{eq:Exp1}
\end{equation}
where we have used the fact that $\mathrm{hor}(\omega^{\alpha})=0$
and $d_{1_{i}}\eta^{\alpha}=\eta_{1_{i}}^{\alpha}$. This expression
(\ref{eq:Exp1}) is independent of the jet-prolongation of the variational
vector field $\eta_{1_{i}}^{\alpha}$ iff 
\[
\rho_{\alpha}^{i}=\partial_{\alpha}^{1_{i}}\mathcal{L}
\]
such that 
\begin{equation}
\mathrm{hor}(j^{1}(\eta)\rfloor\mathrm{d}\mathfrak{C})=\eta^{\alpha}\mathrm{(\partial_{\alpha}\mathcal{L}}-d_{1_{i}}\rho_{\alpha}^{i})\Omega=\eta^{\alpha}\mathrm{\delta_{\alpha}\mathcal{L}}\Omega\label{eq:Exp2}
\end{equation}
holds, where we have used 
\[
\mathrm{hor}(\mathrm{d}\rho_{\alpha}^{i})\wedge\Omega_{i}=d_{1_{i}}(\rho_{\alpha}^{i})\Omega=d_{1_{i}}(\partial_{\alpha}^{1_{i}}\mathcal{L})\Omega.
\]
From (\ref{eq:Exp2}) we see that we have recovered the Euler Lagrange
operator that delivers the domain condition on $\mathcal{S}$, where
$\delta_{\alpha}=\partial_{\alpha}-d_{1_{i}}\partial_{\alpha}^{1_{i}}$
denotes the variational derivative. 

\subsection{Boundary Conditions}

The boundary conditions on $\partial\mathcal{S}$ follow from
\begin{equation}
j^{1}(\eta)\mathfrak{\rfloor C}=\eta^{\alpha}\rho_{\alpha}^{i}\Omega_{i}=\eta^{\alpha}(\partial_{\alpha}^{1_{i}}\mathcal{L})\Omega_{i}.\label{eq:BCPDE1Lag}
\end{equation}
\begin{rem}
\label{rem:BC1}It should be noted that the boundary conditions (\ref{eq:BCPDE1Lag})
on $\partial\mathcal{S}$ can be easily met, if one does not consider
variations on the boundary, i.e. if $\eta$ vanishes. Conversely,
arbitrary variations on $\partial\mathcal{S}$ lead to conditions
for the first-order partial derivatives of $\mathcal{L}$. Furthermore,
in the case of boundary control systems it also possible to extend
the variational problem to include external boundary inputs that equate
with $\partial_{\alpha}^{1_{i}}\mathcal{L}$ such that for arbitrary
variations the boundary conditions are met.
\end{rem}
The functions $\rho_{\alpha}^{j}$ in (\ref{eq:Cartan1odansatz})
have been derived in such a way that in the variational principle
the domain integral becomes independent of derivatives of the variational
field. Furthermore, $\rho_{\alpha}^{j}$ enter the boundary integral
and determine the boundary conditions. As the choice for $\rho_{\alpha}^{j}$
is unique for first-order field theories, so are the expressions in
the boundary integral.

To summarize, (\ref{eq:LagFirstOrderLie}) for $k=1$ is equivalent
to 
\begin{equation}
\int_{\mathcal{S}}j^{2}(\sigma)^{*}(\eta^{\alpha}\delta_{\alpha}\mathcal{L})\Omega+\int_{\partial\mathcal{S}}j^{1}(\sigma)^{*}(\eta^{\alpha}\partial_{\alpha}^{1_{i}}\mathcal{L})\Omega_{i}=0.\label{eq:decompLag1Final}
\end{equation}
\begin{rem}
Alternatively, (\ref{eq:decompLag1Final}) can be derived using the
integration by parts method in a straightforward and unique manner,
as 
\begin{eqnarray*}
\mathrm{L}_{j^{1}(\eta)}(\mathfrak{L}) & = & (\eta^{\alpha}\partial_{\alpha}\mathcal{L}+(d_{1_{i}}\eta^{\alpha})\partial_{\alpha}^{1_{i}}\mathcal{L})\Omega\\
 & = & \eta^{\alpha}\left(\partial_{\alpha}\mathcal{L}-d_{1_{i}}(\partial_{\alpha}^{1_{i}}\mathcal{L})\right)\Omega+d_{1_{i}}(\eta^{\alpha}\,\partial_{\alpha}^{1_{i}}\mathcal{L}\,\Omega)
\end{eqnarray*}
gives the same result as in (\ref{eq:decompLag1Final}) using Stokes'
theorem. Thus, the Cartan form approach indeed qualifies for a coordinate
free version of the integration by parts technique.
\end{rem}

\subsection{Example}

Let us consider the simple model of a Timoshenko beam, see e.g. \cite{Meirovitch},
where the dependent variables $x^{1}$ and $x^{2}$ correspond to
the deflection of a cross section of the beam and the angle of rotation
due to bending. We have that $m=2$ and $r=2$, such that $\mathcal{X}$
possesses the coordinates $(X^{1},X^{2},x^{1},x^{2})$ and the independent
coordinates $(X^{1},X^{2})$ correspond to the time $X^{1}=t$ and
a spatial coordinate $X^{2}=X$. 

The Lagrangian density corresponding to the difference of the kinetic
and potential energy densities is given as $\mathcal{L}\text{\ensuremath{\Omega}\ with \ensuremath{\Omega=\mathrm{d}X^{1}\wedge\mathrm{d}X^{2}}}$
and 
\[
\mathcal{L}=\frac{1}{2}\left(\rho(x_{10}^{1})^{2}+J(x_{10}^{2})^{2}-EI(x_{01}^{2})^{2}-kAG(x_{01}^{1}-x^{2})^{2}\right)
\]
where for simplicity we assume constant beam parameters $\rho,\,J,\,E,\,k,\,A,\,G.$
Using (\ref{eq:decompLag1Final}) the pdes on the domain immediately
follow as 
\begin{equation}
\begin{array}{ccl}
\rho x_{20}^{1} & = & kAG(x_{02}^{1}-x_{01}^{2})\smallskip\\
Jx_{20}^{2} & = & EIx_{02}^{2}+kAG(x_{01}^{1}-x^{2}),
\end{array}\label{eq:TimoPDE}
\end{equation}
as the variational field $\eta$ is arbitrary on $\mathcal{S}$. The
boundary expression follows from the boundary integral in (\ref{eq:decompLag1Final})
and leads to a condition of the form 
\begin{equation}
-\eta^{1}kAG(x_{01}^{1}-x^{2})-\eta^{\text{2}}EIx_{01}^{2}\label{eq:RBtimo}
\end{equation}
which has to vanish at $X=0$ and $X=L$ (the spatial boundary) as
in mechanics no variation on the time boundary takes place. For arbitrary
$\eta^{1}$ and $\eta^{2}$ we derive the classical conditions that
the shear force $kAG(x_{01}^{1}-x^{2})$ and the bending moment $EIx_{01}^{2}$
have to vanish on the spatial boundary.
\begin{rem}
As already pointed out in remark \ref{rem:BC1} it is possible to
extend the variational problem, such that for external boundary inputs
$Q(t)$ and $M(t)$ at $X=0$ and/or $X=L$ the boundary conditions
read as $kAG(x_{01}^{1}-x^{2})=Q$ and $EIx_{01}^{2}=M$. This can
be achieved e.g. by considering an additive boundary integral in the
variational principle such that external boundary inputs are incorporated.
\end{rem}

\section{Second-order Case\label{sec:Second-order-Case}}

Let us consider the case $\mathcal{L}=\mathcal{C}^{\infty}(\mathcal{J}^{2}(\mathcal{X}))$
and consequently the equation (\ref{eq:LagFirstOrderLie}) for $k=2$.
Again, using 
\[
\mathrm{L}_{j^{2}(\eta)}(\mathfrak{L})=j^{2}(\eta)\rfloor\mathrm{d}\mathfrak{L}+\mathrm{d}(j^{2}(\eta)\mathfrak{\rfloor L})
\]
we aim to decompose the variational problem into an expression on
the domain and one on the boundary.

\subsection{The domain conditions}

Similar as in the first-order case we modify the Lagrangian density
by extending it with contact forms of the type 
\[
\omega_{I}^{\alpha}=\mathrm{d}x_{I}^{\alpha}-x_{I+1_{k}}^{\alpha}\mathrm{d}X^{k}\,,\,\,\,0\leq\#I\leq1.
\]
Thus, we make the ansatz 

\begin{eqnarray}
\mathfrak{C} & = & \mathcal{L}\Omega+\rho_{\alpha}^{j}(\mathrm{d}x^{\alpha}-x_{1_{i}}^{\alpha}\mathrm{d}X^{i})\wedge\Omega_{j}+\rho_{\alpha}^{1_{k},j}(\mathrm{d}x_{1_{k}}^{\alpha}-x_{1_{k}+1_{i}}^{\alpha}\mathrm{d}X^{i})\wedge\Omega_{j}\label{eq:Cartan2Ord}
\end{eqnarray}
which reads as
\[
\mathfrak{C}=\mathcal{L}\Omega+\rho_{\alpha}^{j}(\mathrm{d}x^{\alpha}\wedge\Omega_{j}-x_{1_{j}}^{\alpha}\Omega)+\rho_{\alpha}^{1_{k},j}(\mathrm{d}x_{1_{k}}^{\alpha}\wedge\Omega_{j}-x_{1_{k}+1_{j}}^{\alpha}\Omega).
\]
To obtain the domain conditions in an analogous fashion as in the
first-order case we compute $\mathrm{hor}(j^{2}(\eta)\rfloor\mathrm{d}\mathcal{\mathfrak{C}})$
with

\begin{eqnarray*}
\mathrm{d}\mathfrak{C} & = & \mathrm{d}\mathcal{L}\wedge\Omega+\mathrm{d}\rho_{\alpha}^{j}\wedge(\mathrm{d}x^{\alpha}\wedge\Omega_{j}-x_{1_{j}}^{\alpha}\Omega)-\rho_{\alpha}^{j}\mathrm{d}x_{1_{j}}^{\alpha}\wedge\Omega\\
 &  & +\mathrm{d}\rho_{\alpha}^{1_{k},j}\wedge(\mathrm{d}x_{1_{k}}^{\alpha}\wedge\Omega_{j}-x_{1_{k}+1_{j}}^{\alpha}\Omega)-\rho_{\alpha}^{1_{k},j}\mathrm{d}x_{1_{k}+1_{j}}^{\alpha}\wedge\Omega
\end{eqnarray*}
and $\mathrm{hor}(j^{2}(\eta)\rfloor\mathrm{d}\mathcal{\mathfrak{C}})$
can be stated as
\[
\left(\partial_{\alpha}^{I}\mathcal{L}\eta_{I}^{\alpha}-\rho_{\alpha}^{j}\eta_{1_{j}}^{\alpha}-\rho_{\alpha}^{1_{k},j}\eta_{1_{k}+1_{j}}^{\alpha}\right)\Omega-\mathrm{hor}\left(\eta^{\alpha}\mathrm{d}\rho_{\alpha}^{j}+\eta_{1_{k}}^{\alpha}\mathrm{d}\rho_{\alpha}^{1_{k},j}\right)\wedge\Omega_{j}
\]
with $0\leq\#I\leq2$ which simplifies to (\ref{eq:Exp1}) in the
first-order case. Hence, to cancel out the expressions that involve
a jet-prolongation of the variational vector field $\eta_{I}^{\alpha}$,
we have to choose for the case $\#I=2$ 

\begin{equation}
\begin{array}{rrl}
\rho_{\alpha}^{1_{k},j}+\rho_{\alpha}^{1_{j},k} & = & \partial_{\alpha}^{1_{k}+1_{j}}\mathcal{L}\,,\,\,\,k\neq j\\
\rho_{\alpha}^{1_{k},k} & = & \partial_{\alpha}^{1_{k}+1_{k}}\mathcal{L}
\end{array}\label{eq:condCartan2}
\end{equation}
and to guarantee that $\eta_{I}$ vanishes for the case $\#I=1$ we
assign 
\begin{eqnarray}
\rho_{\alpha}^{l} & = & \partial_{\alpha}^{1_{l}}\mathcal{L}-d_{1_{j}}(\rho_{\alpha}^{1_{l},j}).\label{eq:condCartan21}
\end{eqnarray}
Finally, 
\[
\mathrm{hor}(j^{2}(\eta)\rfloor\mathrm{d}\mathfrak{C})=\left(\partial_{\alpha}\mathcal{L}-d_{1_{j}}(\rho_{\alpha}^{j})\right)\eta^{\alpha}\Omega
\]
where we have used the fact that all expressions $\eta_{I}$ for $\#I\neq0$
are canceled out. Using (\ref{eq:condCartan2}) and (\ref{eq:condCartan21})
we obtain 

\[
\partial_{\alpha}\mathcal{L}-d_{1_{l}}(\rho_{\alpha}^{l})=\partial_{\alpha}\mathcal{L}-d_{1_{l}}\left(\partial_{\alpha}^{1_{l}}\mathcal{L}-d_{1_{j}}(\rho_{\alpha}^{1_{l},j})\right)
\]
and consequently 
\[
\mathrm{hor}(j^{2}(\eta)\rfloor\mathrm{d}\mathfrak{C})=\delta_{\alpha}\mathcal{L}\eta^{\alpha}\Omega
\]
with
\begin{equation}
\delta_{\alpha}\mathcal{L}=\sum_{\#I=0}^{2}(-1)^{\#I}d_{I}\partial_{\alpha}^{I}\mathcal{L}.\label{eq:EL2}
\end{equation}
Hence, the domain conditions involve again the variational derivatives
(the components of the Euler-Lagrange operator) (\ref{eq:EL2}) which
is of second-order in this case. The preceding calculation verifies
that despite of the fact that the choice of contact forms in (\ref{eq:condCartan2})
is not unique, it has no influence on the pdes, which is a well-known
fact. But this non uniqueness will have an impact on the boundary
conditions as demonstrated in the following.

\subsection{The boundary conditions}

From (\ref{eq:condCartan2}) it is obvious that there is no unique
choice for the functions $\rho_{\alpha}^{1_{k},j}$ when $k\neq j$
which has the consequence that in the boundary expression, reading
as 
\[
j^{2}(\eta)\mathfrak{\rfloor C}=j^{1}(\eta)\mathfrak{\rfloor C}=\left(\rho_{\alpha}^{i}\eta^{\alpha}+\rho_{\alpha}^{1_{k},i}\eta_{1_{k}}^{\alpha}\right)\Omega_{i},
\]
a different choice for $\rho_{\alpha}^{1_{k},i}$ leads to different
boundary conditions. To overcome this problem we pose an extra condition
on the set of possible boundary terms in order to extract the physical
meaningful boundary terms. To this end, we assume that the boundary
of the submanifold $\mathcal{S}$ which is denoted as $\partial\mathcal{S}$
is parameterized as $X^{r}=const.$ such that the coordinates for
$\mathcal{D}$ are $(X_{\partial}^{i},\,X^{r})$ with $X_{\partial}^{i}=(X^{1},\ldots,X^{r-1})$
which is possible at least locally by a change of the independent
coordinates, such that in these coordinates adapted to the boundary
the boundary volume form $\Omega_{\partial}=\partial_{r}\rfloor\Omega$
can be introduced. It is now remarkable that for a section $\sigma_{\partial}$
on the boundary, i.e. $x^{\alpha}=\sigma_{\partial}^{\alpha}(X_{\partial}^{i})$,
we can compute all the first-order jet-variables apart from the $X^{r}$
direction, which represents an additional degree of freedom on the
boundary. 
\begin{rem}
Let us consider the inclusion mapping $\iota:\partial\mathcal{D}\rightarrow\mathcal{D}$
in adapted coordinates. Then $\iota^{*}(\mathcal{X})$ and $\iota^{*}(\mathcal{J}^{1}(\mathcal{X}))$
are equipped with induced coordinates $(X_{\partial}^{i},x^{\alpha})$
and $(X_{\partial}^{i},x^{\alpha},x_{1_{l}}^{\alpha})\,,\,l=1,\ldots,r$,
respectively. It is worth noting that $\iota^{*}(\mathcal{J}^{1}(\mathcal{X}))$
does not qualify as a jet-manifold as we have too few independent
coordinates, but the coordinate $x_{1_{r}}^{\alpha}$ will play an
essential role in deriving the boundary conditions. This construction
has been used also in \cite{EnnsPhD} to introduce the so-called extended
Cartan form on the boundary. We will pursue a different target, namely
we will exploit this boundary structure to give conditions for the
choice of the under determined functions $\rho_{\alpha}^{1_{k},i}$.
\end{rem}
The boundary conditions are derived by means of 
\begin{equation}
j^{1}(\eta)\mathfrak{\rfloor C}=\left(\rho_{\alpha}^{r}\eta^{\alpha}+\rho_{\alpha}^{1_{k},r}\eta_{1_{k}}^{\alpha}\right)\Omega_{\partial}=\left(\rho_{\alpha}^{r}\eta^{\alpha}+\rho_{\alpha}^{1_{1},r}\eta_{1_{1}}^{\alpha}+\ldots+\rho_{\alpha}^{1_{r},r}\eta_{1_{r}}^{\alpha}\right)\Omega_{\partial}\label{eq:BoundCartan}
\end{equation}
where the summation is carried out over $k$ as $r$ is fixed in the
adapted coordinates. On the boundary, independent variations with
respect to $x^{\alpha}$ and $x_{1_{r}}^{\alpha}$ can be carried
out, as discussed above, such that we are led to set 
\[
\rho_{\alpha}^{1_{k},r}=0\,,\,\,\,\,k=1,\ldots,r-1
\]
such that
\begin{equation}
j^{1}(\eta)\mathfrak{\rfloor C}=\left(\rho_{\alpha}^{r}\eta^{\alpha}+\rho_{\alpha}^{1_{r},r}\eta_{1_{r}}^{\alpha}\right)\Omega_{\partial}\label{eq:BoundCartanPartial}
\end{equation}
is met. 
\begin{rem}
From a variational point of view $\rho_{\alpha}^{1_{k},r}=0,\,k\neq r$
guarantees that only independent variations arise in the boundary
integral, which can be interpreted as an additional postulate which
selects the appropriate Cartan form out of a set of possible Cartan
forms that all deliver the same domain conditions. It should be noted
that this choice always leads to the minimal possible boundary conditions.
The coordinate representation is well justified in adapted coordinates
only, i.e. when the boundary is parameterized as $X^{r}=const.$
\end{rem}
This means that in the relations (\ref{eq:condCartan2}) we have to
choose $\rho_{\alpha}^{1_{r},k}=\partial_{\alpha}^{1_{r}+1_{k}}\mathcal{L}$.
From relation (\ref{eq:condCartan21}) we observe that then the coefficient
of $\eta^{\alpha}$ namely $\rho_{\alpha}^{r}$ which is 
\begin{equation}
\rho_{\alpha}^{r}=\partial_{\alpha}^{1_{r}}\mathcal{L}-d_{1_{k}}(\rho_{\alpha}^{1_{r},k}).\label{eq:rhoboud}
\end{equation}
involves exactly the terms $\rho_{\alpha}^{1_{r},k}$. Hence we have
exploited the degree of freedom regarding the functions $\rho_{\alpha}^{1_{k},i}$
in (\ref{eq:condCartan2}) to generate the admissible boundary expressions. 
\begin{rem}
It should be noted that in \cite{EnnsPhD} a different strategy is
suggested which requires an additional integration by parts or in
terms of \cite{EnnsPhD} an extended Cartan form which contains an
additive total differential which corrects the wrong boundary terms.
E.g. an expression involving $\eta_{1_{k}}^{\alpha}$ for $1_{k}\neq1_{r}$
can be again integrated by parts to produce an admissible $\eta^{\alpha}$
term. Both strategies give the same boundary conditions as the additional
integration by parts exactly delivers the second term in the right
hand side of (\ref{eq:rhoboud}).
\end{rem}

\section{Examples for second-order mechanical systems}

In this section we demonstrate our approach by means of two examples
- the Euler-Bernoulli beam and the Kirchhoff plate. 

\subsection{Euler-Bernoulli Beam}

We have two independent variables $X^{1}=t$ and $X^{2}=X$, hence
$\mathrm{dim}(\mathcal{D})=2$, where $t$ is the time coordinate
and $X$ the spatial variable (we set $\mathcal{S}=\mathcal{D}$ for
simplicity). The dependent variable is the deflection $x^{1}$, such
that $\mathcal{X}$ is equipped with the coordinates $(X^{1},X^{2},x^{1})$
and $\Omega=\mathrm{d}X^{1}\wedge\mathrm{d}X^{2}=\mathrm{d}t\wedge\mathrm{d}X$
together with the Lagrangian (kinetic minus potential energy)
\[
\mathcal{L}=\frac{1}{2}\rho A(x_{10}^{1})^{2}-\frac{1}{2}EI(x_{02}^{1})^{2}
\]
where again for simplicity we consider constant physical parameters
$\rho,\,A,\,E,\,I$. The pdes follow from the Euler-Lagrange operator
(\ref{eq:EL2}), which in this case simplifies to $\delta_{1}=-d_{10}\partial_{1}^{10}+d_{02}\partial_{1}^{02}$
such that we obtain the well-known pde
\[
\rho Ax_{20}^{1}+EIx_{04}^{1}=0.
\]

\subsubsection{Boundary Conditions}

To obtain the boundary conditions, we consider the relation (\ref{eq:BoundCartan})
with $\Omega_{\partial}=\Omega_{r}=\partial_{2}\rfloor\Omega=-\mathrm{d}t$
since we are interested in the spatial boundary (no variation on the
time boundary) and derive
\[
\left(\rho_{1}^{2}\eta^{1}+\rho_{1}^{10,2}\eta_{10}^{1}+\rho_{1}^{01,2}\eta_{01}^{1}\right)\Omega_{\partial}=0.
\]
As no mixed partial derivatives appear in $\mathcal{L}$, $\rho_{1}^{10,2}=0$
is met automatically and we have 
\begin{eqnarray*}
\rho_{1}^{01,2} & = & \partial_{1}^{02}\mathcal{L}=-EIx_{02}^{1}\\
\rho_{1}^{2} & = & -d_{01}\rho_{1}^{01,2}=EIx_{03}^{1}
\end{eqnarray*}
such that we obtain the conditions on the spatial boundary as 
\[
\eta^{1}EIx_{03}^{1}=0\,,\,\,\,\eta_{01}^{1}EIx_{02}^{1}=0.
\]
Here, $EIx_{03}^{1}$ corresponds to the shear force and $EIx_{02}^{1}$
to the bending moment and $\eta^{1}$ and $\eta_{01}^{1}$ take their
values according to the support on the boundary.
\begin{rem}
\label{remBCEulerB}Let us consider a cantilever beam that is built-in
at $X=0$ (clamped end) and free at $X=L$ (free end). Then at $X=0$
the imposed geometric boundary conditions imply $\eta^{1}=\eta_{01}^{1}=0$
whereas at $X=L$ we have the dynamic boundary conditions $EIx_{02}^{1}=EIx_{03}^{1}=0$.
Furthermore, it is again possible to apply external forces and torques
at $X=L$ which can be included in the variational principle. 
\end{rem}
As no mixed partial derivatives appear in this example, the Cartan
form approach is unambiguous although the example belongs to second-order
theory. This will be different for the following example.

\subsection{Kirchhoff Plate\label{subsec:Kirchhoff-Plate}}

Let us consider the case of a (rectangular) Kirchhoff plate with three
independent variables $X^{1}=t,\,X^{2}=X$ and $X^{3}=Y$, such that
we have two spatial variables and the time. Hence, $\mathcal{X}$
is equipped with the coordinates $(X^{1},X^{2},X^{3},x^{1})$ where
$x^{1}$ denotes the deflection. The Lagrangian is given as $\mathcal{L}=\mathcal{T}-\mathcal{V}$
with $2\mathcal{T}=(x_{100}^{1})^{2}$ and 
\[
2\mathcal{V}=(x_{020}^{1})^{2}+(x_{002}^{1})^{2}+2\nu x_{020}^{1}x_{002}^{1}+2(1-\nu)(x_{011}^{1})^{2}
\]
where we have set the physical parameters to one apart from the Poisson's
ratio $\nu$. Again the pdes follow from the Euler-Lagrange operator
(\ref{eq:EL2}), with 
\[
\delta_{1}=-d_{100}\partial_{1}^{100}+d_{020}\partial_{1}^{020}+d_{002}\partial_{1}^{002}+d_{011}\partial_{1}^{011}
\]
such that we obtain
\begin{equation}
x_{200}^{1}+x_{040}^{1}+x_{004}^{1}+2x_{022}^{1}=0.\label{eq:pdeKirch}
\end{equation}

\subsubsection{The boundary conditions}

To obtain the boundary conditions on the boundary $X^{3}=Y=const.$
we consider
\begin{equation}
\left(\rho_{1}^{3}\eta^{1}+\rho_{1}^{010,3}\eta_{010}^{1}+\rho_{1}^{001,3}\eta_{001}^{\alpha}\right)\Omega_{\partial}\label{eq:Kirchbound1}
\end{equation}
where we already have used $\rho_{1}^{100,3}=0$ since no mixed time
and spatial derivatives appear. The non-zero components for $\rho_{1}$
of regarding the jet-order $\#I=2$ in $\mathcal{L}$ are
\[
\rho_{1}^{001,3}=\partial_{1}^{002}\mathcal{L}\,,\,\,\,\rho_{1}^{010,2}=\partial_{1}^{020}\mathcal{L}
\]
as well as
\[
\rho_{1}^{010,3}+\rho_{1}^{001,2}=\partial_{1}^{011}\mathcal{L}.
\]
It is remarkable that for the mixed spatial derivative we get two
possible contact forms. For $\#I=1$ we obtain from $\rho_{\alpha}^{l}=\partial_{\alpha}^{1_{l}}\mathcal{L}-d_{1_{j}}(\rho_{\alpha}^{1_{l},j})$
the coordinate expressions
\begin{eqnarray*}
\rho_{1}^{3} & = & \partial_{1}^{001}\mathcal{L}-d_{001}\rho_{1}^{001,3}-d_{010}\rho_{1}^{001,2}\\
\rho_{1}^{2} & = & \partial_{1}^{010}\mathcal{L}-d_{010}\rho_{1}^{010,2}-d_{001}\rho_{1}^{010,3}
\end{eqnarray*}
and $\rho_{1}^{1}=\partial_{1}^{100}\mathcal{L}$. Since $\eta_{010}^{1}$
is not admissible on the boundary of interest we can choose $\rho_{1}^{010,3}=0$
and consequently $\rho_{1}^{001,2}=\partial_{1}^{011}\mathcal{L}$
follows. Hence,
\begin{equation}
\rho_{1}^{3}=-d_{001}\rho_{1}^{001,3}-d_{010}\rho_{1}^{001,2}=(x_{003}^{1}+(2-\nu)x_{021}^{1})\label{eq:Kirch000term}
\end{equation}
and therefore on the spatial boundary $X^{3}=Y=const.$ we obtain
(no variation on the time boundary)
\begin{eqnarray}
\eta^{1}(x_{003}^{1}+(2-\nu)x_{021}^{1}) & = & 0\label{eq:boundk1}\\
-\eta_{001}^{1}(x_{002}^{1}+\nu x_{020}^{1}) & = & 0\label{eq:boundk2}
\end{eqnarray}
where (\ref{eq:boundk1}) is a condition for the shear force and (\ref{eq:boundk2})
for the bending moment. As a consequence, we obtain two boundary conditions,
whereas a wrong choice for $\rho_{1}$ (corresponding to a wrong integration
by parts) would deliver three (compare with (\ref{eq:Kirchbound1})).
To construct the boundary conditions for the boundary $X^{2}=X=const.$
one can easily adopt the procedure or re-parameterize the independent
coordinates. Furthermore, it should be noted that in an analogous
fashion as in remark \ref{remBCEulerB} the introduction of imposed
geometric and dynamic boundary conditions can be accomplished.

\section{Boundary Ports }

To motivate for the idea of boundary ports, let et us consider a Hamiltonian
of first-order based on an underlying bundle structure $\mathcal{X}\rightarrow\mathcal{D}$.
Hence, we have $\mathfrak{H}=\mathcal{H}\Omega$ with $\mathcal{H}\in C^{\infty}(\mathcal{J}^{1}(\mathcal{X}))$.
Depending on a chosen section $s:\mathcal{D}\rightarrow\mathcal{X}$
we are led to compute $H(s)=\int_{\mathcal{D}}j^{1}(s)^{*}(\mathfrak{H})$.
Let us consider now an evolutionary vector-field $v=v^{\alpha}\partial_{\alpha}$,
corresponding to the pdes
\begin{equation}
\dot{x}^{\alpha}=v^{\alpha}\,,\,\,\,\,v^{\alpha}\in C^{\infty}(\mathcal{J}^{2}(\mathcal{X}))\label{eq:EvolPDE}
\end{equation}
and we are interested in the change of $H$ along the section $s$,
where $s$ is a formal solution of the pdes parameterized in $t$,
i.e.
\[
\partial_{t}s(X)=v(X,\partial_{J}s(X))\,,\,\,\,\#J\leq2
\]
together with appropriate boundary conditions, where we have suppressed
the explicit dependence of $s$ with respect to the evolution parameter
and well-posedness of the problem is assumed.

The idea of ports in a Hamiltonian setting can be visualized by computing
a power balance relation by means of the Hamiltonian density which
can be interpreted as an energy density in many applications. Thus
we evaluate 
\begin{equation}
\dot{H}=\partial_{t}H(s)=\int_{\mathcal{S}}j^{2}(s){}^{*}(\mathrm{L}_{j^{1}(v)}\mathfrak{H)}\label{eq:LagFirstOrderLie-1}
\end{equation}
and regarding the evaluation of the expression $\mathrm{L}_{j^{1}(v)}\mathfrak{H}$
appearing in (\ref{eq:LagFirstOrderLie-1}) the same reasoning as
in section \ref{sec:First-order-Case} can be applied and 
\begin{equation}
\partial_{t}H(s)=\int_{\mathcal{S}}j^{2}(s)^{*}(v^{\alpha}(\delta_{\alpha}\mathcal{H})\Omega)+\int_{\partial\mathcal{S}}j^{1}(s)^{*}(v^{\alpha}(\partial_{\alpha}^{1_{i}}\mathcal{H})\Omega_{i})\label{eq:Decomp1}
\end{equation}
follows which is the counterpart to (\ref{eq:decompLag1Final}).
\begin{rem}
It should be noted that instead of a variational vector field $\eta$
that performs the variation, now an evolutionary vector field that
corresponds to a dynamical systems is applied. But form a mathematical
point of view the calculations remain valid when $\mathcal{L}$ is
replaced by $\mathcal{H}$ and $\eta$ is replaced by $v$. 
\end{rem}
Thus, we have obtained a power balance relation when $\mathfrak{H}=\mathcal{H}\Omega$
corresponds to an energy density. Obviously, in this power balance
relation we have a contribution in the domain and one on/over the
boundary. We can restate (\ref{eq:Decomp1}) to obtain 
\begin{equation}
\dot{H}=\int_{\mathcal{S}}j^{2}(s)^{*}(v\rfloor\delta\mathfrak{H})+\int_{\partial\mathcal{S}}j^{1}(s)^{*}(v\rfloor\delta^{\partial}\mathfrak{H})\label{eq:PowerBal1st}
\end{equation}
with $\delta\mathfrak{H}=(\delta_{\alpha}\mathcal{H})\mathrm{d}x^{\alpha}\wedge\Omega$
and $\delta^{\partial}\mathfrak{H}=(\partial_{\alpha}^{1_{i}}\mathcal{H})\mathrm{d}x^{\alpha}\wedge\Omega_{i}$.
Hence, choosing
\begin{equation}
v=(\mathcal{J}-\mathcal{R})(\delta\mathfrak{H})+u\rfloor\mathcal{G}\label{eq:Hamv1}
\end{equation}
where $\mathcal{J}$ is a skew- symmetric map and $\mathcal{R}$ is
positive semi-definite, we can express the impact on the domain as
\begin{equation}
\int_{\mathcal{S}}j^{2}(s)^{*}(-\mathcal{R}(\delta\mathfrak{H})\rfloor\delta\mathfrak{H}+y\rfloor u)=\int_{\mathcal{S}}j^{2}(s)^{*}(-(\delta_{\alpha}\mathcal{H})\mathcal{R}^{\alpha\beta}(\delta_{\beta}\mathcal{H})+y_{\xi}u^{\xi})\Omega\label{eq:PowerBal1sta}
\end{equation}
with $y=\mathcal{G}^{*}(\delta\mathfrak{H})$ and the input map $\mathcal{G}$
together with its adjoint $\mathcal{G}^{*}.$ The boundary impact
(depending on the boundary conditions of the pdes) is determined by
the expression $v\rfloor\delta^{\partial}\mathfrak{H}$.
\begin{cor}
Given evolutionary pdes of the form (\ref{eq:EvolPDE}) with the generalized
vector field $v$ as in (\ref{eq:Hamv1}), the formal change of the
first-order Hamiltonian $H$ along solutions of (\ref{eq:EvolPDE})
can be decomposed in a dissipative term in the domain $-\mathcal{R}(\delta\mathfrak{H})\rfloor\delta\mathfrak{H}$,
in a power port in the domain $y\rfloor u$ and in a power port on
the boundary $v\rfloor\delta^{\partial}\mathfrak{H}.$
\end{cor}
The representation (\ref{eq:Hamv1}) together with the collocated
output $y=\mathcal{G}^{*}(\delta\mathfrak{H})$ can be seen as a generalization
of the classical port-Hamiltonian formulation of ODEs, see e.g. \cite{vanderSchaft},
to the infinite-dimensional case based on the underlying jet-bundle
structure and has been proposed and discussed in \cite{SchlacherHam2008,SchoeberlTAC2013,SchoeberlAut2014PDE}.
The power balance relation (\ref{eq:PowerBal1sta}) is valid when
the maps $\mathcal{J},\,\mathcal{R}$ and $\mathcal{G}$ do not involve
differential operators. However, also this differential operator case
can be discussed in the proposed setting, see e.g. \cite{SchoeberlAut2014PDE},
but then additional boundary expressions have to be taken into account
due to the differential operators. The first-order case can be also
handled by using the classical integration by parts technique, as
it has been done in the cited references. For second-order theories
similar problems as in the Lagrangian setting (calculus of variations)
appear and we will show that the proposed Cartan form approach can
also be exploited in the port-Hamiltonian setting for the determination
of the boundary ports in a similar way as it has been used in the
Lagrangian scenario to determine the boundary conditions. 

\subsection{Domain and Boundary Impact for the second-order case}

Let us consider a Hamiltonian density $\mathfrak{H}=\mathcal{H}\Omega$
with $\mathcal{H}\in C^{\infty}(\mathcal{J}^{2}(\mathcal{X}))$ and
an evolutionary vector field $v$. Then, to evaluate 
\[
\dot{H}=\partial_{t}H(s)=\int_{\mathcal{S}}j^{4}(s)^{*}(\mathrm{L}_{j^{2}(v)}\mathfrak{H)}
\]
we follow the same consideration as in section \ref{sec:Second-order-Case}
by augmenting $\mathcal{H}\Omega$ with appropriate contact forms.
This immediately gives us the desired decomposition of $\dot{H}$
into a domain expression $\dot{H}_{\mathcal{S}}$ and a boundary expression
$\dot{H}_{\partial\mathcal{\mathcal{S}}}$ with 

\[
\dot{H}_{\mathcal{S}}=\int_{\mathcal{S}}j^{4}(s)^{*}(v\rfloor\delta\mathfrak{H})
\]
where $\delta$ involves the Euler-Lagrange operator of second-order,
i.e. 
\begin{equation}
\delta_{\alpha}\mathcal{H}=\sum_{\#I=0}^{2}(-1)^{\#I}d_{I}\partial_{\alpha}^{I}\mathcal{H}.\label{eq:EulerLagforH}
\end{equation}
Thus, again the consideration of (\ref{eq:Hamv1}), where $\delta\mathfrak{H}$
involves (\ref{eq:EulerLagforH}) in this case, gives an analogous
power balance relation. The crucial point is again the boundary impact
$\dot{H}_{\partial\mathcal{\mathcal{S}}}$ and the key to solve the
problem is to consider adapted coordinates with respect to the boundary
such that 
\begin{equation}
\dot{H}_{\partial\mathcal{\mathcal{S}}}=\int_{\partial\mathcal{S}}j^{3}(s)^{*}\left(\rho_{\alpha}^{r}v^{\alpha}+\rho_{\alpha}^{1_{r},r}v_{1_{r}}^{\alpha}\right)\Omega_{\partial}\label{eq:BunfHamPartial}
\end{equation}
can be stated (which is the counterpart to (\ref{eq:BoundCartanPartial})),
where the determination of the function $\rho_{\alpha}^{r}$ and $\rho_{\alpha}^{1_{r},r}$
follows from (\ref{eq:condCartan2}) and (\ref{eq:condCartan21})
by substituting $\mathcal{L}$ with $\mathcal{H}.$

\subsection{Example}

The port-Hamiltonian representation of the Kirchhoff plate can be
accomplished by choosing a manifold $\mathcal{X}$ equipped with the
coordinates $(X^{1},X^{2},x^{1},p_{1})$ where $x^{1}$ again denotes
the deflection and $x^{2}=p_{1}$ the linear momentum It should be
noted that the number of independent variables is two, in contrast
to section \ref{subsec:Kirchhoff-Plate}, as the time is serving as
an evolution parameter. The Hamiltonian is given as $\mathcal{H}=\mathcal{T}+\mathcal{V}$
with $2\mathcal{T}=(p_{1})^{2}$ and 
\[
2\mathcal{V}=(x_{20}^{1})^{2}+(x_{02}^{1})^{2}+2\nu x_{20}^{1}x_{02}^{1}+2(1-\nu)(x_{11}^{1})^{2}.
\]
The pdes read as
\begin{equation}
\left[\begin{array}{c}
\dot{x}^{1}\\
\dot{p}_{1}
\end{array}\right]=\left[\begin{array}{cc}
0 & 1\\
-1 & 0
\end{array}\right]\left[\begin{array}{c}
\delta_{1}\mathcal{H}\\
\delta^{1}\mathcal{H}
\end{array}\right]\,,\,\,\,\,\mathcal{J}=\left[\begin{array}{cc}
0 & 1\\
-1 & 0
\end{array}\right]\label{eq:KirchhoffHamEx}
\end{equation}
with 
\[
\delta^{1}=\frac{\partial}{\partial p_{1}}\,,\,\,\,\delta_{1}=\partial_{1}-d_{1_{i}}\partial_{1}^{1_{i}}+d_{J}\partial_{1}^{J}\,,\,\,\,\,\#J=2
\]
which is
\[
\dot{x}^{1}=p_{1}\,,\,\,\,\,\dot{p}_{1}=-(x_{40}^{1}+x_{04}^{1}+2x_{22}^{1}).
\]
The Hamiltonian includes second-order jet variables, and the application
of the variational derivative implies then forth-order pdes in our
evolutionary setting.
\begin{rem}
The pdes correspond to (\ref{eq:pdeKirch}) where the second time
derivative, i.e. $x_{200}$, is written as a first-order derivative
with respect to the momentum variable $p_{1}$ to obtain an evolutionary
representation in the Hamiltonian setting. Furthermore, it is easily
checked that (\ref{eq:KirchhoffHamEx}) is of the form $v=(\mathcal{J}-\mathcal{R})(\delta\mathfrak{H})+u\rfloor\mathcal{G}$
with $\mathcal{R}=0$ and $\mathcal{G}=0$ and with the second-order
variational derivatives $\delta^{1}$ and $\delta_{1}$ acting on
$\mathcal{H}$, where $\delta^{1}$ degenerates to a partial derivative
as no jet-variables with respect to $p_{1}$ are involved.
\end{rem}
Next we wish to analyze the power balance relation in order to obtain
the corresponding ports. As we use a skew-symmetric map $\mathcal{J}$
we can conclude that $\dot{H}_{\mathcal{D}}$ in $\dot{H}=\dot{H}_{\mathcal{D}}+\dot{H}_{\partial\mathcal{D}}$
is equal to zero and that only a boundary impact has to be considered.

\subsubsection{The boundary impact}

From (\ref{eq:BunfHamPartial}) we obtain on the boundary $X^{2}=const$.
the expression

\[
\dot{H}_{\partial\mathcal{D}}=-\int(v^{1}\rho_{1}^{2}+v_{01}^{1}\rho_{1}^{01,2})\mathrm{d}X^{1}=\int\left(\dot{x}^{1}F_{S}+\dot{x}_{01}^{1}M_{B}\right)\mathrm{d}X^{1}
\]
with the shear force $F_{S}=(x_{03}^{1}+(2-\nu)x_{21}^{1})$ appearing
in the power port $\dot{x}^{1}F_{S}$ and the bending moment $M_{B}=-(x_{02}^{1}+\nu x_{20}^{1})$
in $\dot{x}_{01}^{1}M_{B}$. The collocated quantities are $\dot{x}^{1}=v^{1}$
(velocity) and $\dot{x}_{01}^{1}=v_{01}^{1}$ (angular velocity).
Hence, to obtain the physically correct ports, it is vitally important
to perform the integration by parts properly, which in our language
amounts for a correct selection of the coefficients of the Cartan-form
in accordance with the parameterization of the boundary.

\section{Conclusion}

In this contribution we have provided a strategy to compute the boundary
conditions for second-order field theories in a variational formulation
in an algorithmic fashion suitable for the computer algebra implementation.
Furthermore, for the class of port-Hamiltonian systems formulated
on jet-bundles, the same reasoning regarding a minimal number of boundary
ports can be carried out, as the mathematical problem of computing
the variation of the Lagrangian density is equivalent to the evaluation
of the change of the Hamiltonian (the total energy) along solutions
of an evolutionary set of pdes in a port-Hamiltonian setting.

\bibliographystyle{plain}
\bibliography{maxibib}

\end{document}